\documentclass[12pt]{article}
\title{\bf Absolute Primes}
\author{\bf Arkadii Slinko}

\date{Department of Mathematics\\
The University of Auckland\\
Private Bag 92019\\
Auckland 1142, New Zealand\\
a.slinko@auckland.ac.nz }

\begin{document}
\maketitle

For a long time primes have attracted the attention of mathematicians, especially those
primes that possess some sort of symmetry. The mysterious and inconceivable repunits
\[
A_n={111\dots 1}_{(n)},
\]
whose decimal representation contains only units, form an important
class of them\footnote{By a subscript in brackets we will indicate the number of digits
in the decimal representation of the integer.}.  For a repunit $A_n$ to be prime, the
number $n$ of digits in its decimal representation must be also prime. But this
condition is far from being sufficient: for instance, $A_3=111=3\cdot
37$ and $A_5=11111=41\cdot 271$. Some repunits are nonetheless prime:
$A_2$, $A_{19}$, $A_{23}$, $A_{317}$ and possibly $A_{1031}$, give the
examples.  The question of primeness of repunits was discussed by M.
Gardner [1] and later in [2-4]. It is completely unclear whether the
number of prime repunits is finite or infinite.

The prime repunits are examples of integers which are prime and remain prime
after an arbitrary permutation of their decimal digits. Integers with
this property are called either {\em permutable primes} according to
H.-E. Richert  [5], who introduced them some 40 years ago, or {\em
absolute primes} according to T.N. Bhagava and P.H. Doyle [6], and
A.W.Johnson [7]. The intent of this note is to give a short proof, which
does not require significant number crunching, of all known facts referring to
absolute primes different from repunits.    
\bigskip

Analyzing the table of primes which are less than $10^3$,  we find $21$
absolute primes different form the repunit $11$:
\[
2,3,5,7,13,17,31,37,71,73,79,97,113,131,199,311,337,373,733,919,991.
\]  

The first observation is easy to get:
\medskip

{\bf Lemma 1:} A multidigit absolute prime contains in its decimal
representations only the four digits $1$, $3$, $7$, $9$.

{\bf Proof:} If any of the digits $0$, $2$, $4$, $5$, $6$,
$8$  appear in the representation of an integer, then by shifting this digit to
the units place we get a multiple of $2$ or a multiple of $5$.
\medskip

Now we can confine the area of the search, and this helps us to proceed with the
following deliberations.       
\medskip

{\bf Lemma 2:} An absolute prime does not contain in its decimal representation all of
the digits $1$, $3$, $7$, $9$ simultaneously. 

{\bf Proof:} Let $N$ be a number with all of the digits $1$, $3$, $7$, $9$ in its
decimal representation. Let us shift these four digits to the rightmost four places, to
obtain an integer
\[
N_0=\overline{c_1\dots c_{n-4}7931}=L\cdot 10^5+7931,
\]
where the notation $\overline{a_1\dots a_{n}}$ is used to denote the number
$a_110^{n-1}+a_210^{n-2}+\dots +a_{n-1}10+a_n$, which decimal representation consists of
digits $a_1,\dots ,a_n$. The integers $K_0=7931$, $K_1=1793$, $K_2=9137$, $K_3=7913$,
$K_4=7193$, $K_5=1937$, $K_6=7139$ have different remainders on dividing by $7$ since
$K_i\equiv i\ \mbox{(mod 7)}$. The seven integers $N_i=L\cdot 10^5+K_i$ for
$i=0,1,2,3,4,5,6$ also have different remainders on dividing by $7$. Therefore one of
them is a multiple of
$7$. Since these integers can be obtained from $N$ by a permutation of digits, $N$ is not
an absolute prime.
\medskip

{\bf Lemma 3:} No absolute prime contains in its decimal representation three digits $a$
and two digits $b$ simultaneously, provided $a\ne b$.

{\bf Proof:} Suppose that an integer $N$ contains digits $a,a,a,b,b$ in its decimal
representation. By a permutation of digits of $N$, we can obtain integers
\[
N_{i,j}=\overline{c_1\dots c_{n-5}aaaaa}+(b-a)(10^i+10^j),
\]
where $4\ge i>j\ge 0$. Since the integers $10^4+10^1$, $10^3+10^2$, $10^3+10^1$,
$10^2+10^0$, $10^1+10^0$, $10^4+10^0$, $10^4+10^2$ yield different remainders on
dividing by $7$, which are respectively $0,1,2,3,4,5,6$, so do the integers
$(b-a)(10^i+10^j)$, when $4\ge i>j\ge 0$. Therefore among the integers $N_{i,j}$ there
exists one that is a multiple of $7$. 
\medskip

Using these two lemmas we are able by direct calculation (preferably on a computer) to
prove that no
$n$-digit absolute primes exist with $n=4,5,6$. For example, if $n=6$, we have to
check all the numbers
\[
\overline{aaaaab},\quad \overline{aaaabc},\quad \overline{aabbcc},
\]
where $a,b,c\in\{1,3,7,9\}$.

{\bf Lemma 4:} If  $N=\overline{c_1\dots c_{n-6}aaaaab}$ is an absolute prime, then $K=
\overline{c_1\dots c_{n-6}}$ is divisible by $7$.

{\bf Proof:} By permutations of the last $6$ digits of $N$ we can obtain the integers 
\[
N_i=K\cdot 10^6+a\cdot A_6+(b-a)10^i
\]
for $0\le i\le 5$. Since $b-a$ is even and the powers $10^i$, $0\le i\le 5$, have
different nonzero remainders on dividing by $7$:
\[
10^0\equiv 1,\quad 10^1\equiv 3,\quad 10^2\equiv 2,\quad 10^3\equiv
6,\quad 10^4\equiv  4,\quad 10^5\equiv 5,  
\]
the integers $(b-a)10^i$ have the same property. If the integer $K\cdot 10^6+a\cdot A_6$
had a nonzero remainder on dividing by $7$, we would find an integer $(b-a)10^i$, which
has the opposite remainder, and get $N_i$ which is divisible by $7$. Since this is
impossible, the number $K\cdot 10^6+a\cdot A_6$ is a multiple of $7$. Moreover, as
$A_6=10^0+10^1+10^2+10^3+10^4+10^5\equiv 1+3+2+6+4+5=21\equiv 0\ \mbox{(mod 7)}$, we
conclude that $K\cdot 10^6$ and hence $K$, is divisible by $7$.
\medskip

{\bf Theorem 1:} Every multidigit absolute prime integer $N$ is either a repunit, or
can be obtained by a permutations of digits from the integer
\[
B_n(a,b)=\overline{aaa\dots
ab}_{(n)}=a\cdot A_n+(b-a),                                                        
\]
where $a$ and $b$ are different digits from the set $\{1,3,7,9\}$.

{\bf Proof:} Let $n$ be the number of digits of $N$.  We  can  suppose 
that $n>6$. By the first three  lemmas $N$ is written by the digits $1,3,7,9$ only but it
does  not  contain  in  its  decimal representation all of the digits $1,3,7,9$, and
it can contain three such digits only if $N$  is  a permutation of digits of the  number 
$\overline{aaa\dots abc}_{(n)}$.  Let us  show  that this is impossible. Since $N$ is an
absolute prime, the integers
\[
N_1 =\overline{a\dots acaaaaab}_{(n)},\qquad N_2 =\overline{a\dots
abaaaaac}_{(n)},                 
\]
are also absolute primes, and by Lemma 4  the  integers  $\overline{a\dots
ac}_{(n-6)}$  and $\overline{a\dots ab}_{(n-6)}$  are both divisible by $7$. Thus their 
difference,  whose  absolute value is $|b-c|$, is also divisible by $7$, which is
impossible.

Therefore either $N$ is a repunit or else it is written by two digits only. 
In the latter case we need Lemma 3 again, to secure that one  digit 
appears only once. 
\medskip

The prime number $7$ played a significant role in 
the preceding considerations. But other useful primes also exist 
and we are going to find some of them. Note that the  property of $7$ most useful 
for us was the fact that the powers $10^i$, $0<i<6$,  had 
different nonzero remainders on  dividing  by  $7$.  In  general, by 
Fermat's Little Theorem  for every  prime  $p>5$,  we  have 
$10^{p-1}\equiv 1\ \mbox{(mod $p$)}$. 

Let $h(p)$ be the least  possible  positive  integer  such  that 
$10^{h(p)}\equiv 1\ \mbox{(mod $p$)}$. It is obvious that $h(p)$ is a divisor of $p-1$ and 
that 
$10^{q}\equiv 1\ \mbox{(mod $p$)}$ implies that $q$ is divisible by $h(p)$. It is also easy to see
that the powers $10^j$, $0<j<p-1$, have different nonzero remainders on dividing  by $p$ 
as soon as $h(p)=p-1$. When this is the case, $10$ is  said  to  be  a 
{\bf primitive root} modulo $p$.
 
Note that $10$ is a primitive root modulo primes $17$, 
$19$, $23$, $29$ (the reader again may write a computer program to check that), but
$10$  is  not  a  primitive  root  modulo $13$  since  $10^6\equiv 1\ \mbox{(mod $13$)}$.
\medskip

{\bf Lemma 5:} Let $A_n$  be a repunit and $p>3$ is a prime. Then $A_n\equiv 0\ \mbox{(mod $p$)}$ 
if and only if $n\equiv 0\ \mbox{(mod $h(p)$)}$.

{\bf Proof:} As $10^n  =9\cdot A_n + 1$, we have $A_n\equiv 0\ \mbox{(mod $p$)}$ if and only if
$10^n \equiv 1\ \mbox{(mod $p$)}$ and this is equivalent to $n\equiv 0\ \mbox{(mod $h(p)$)}$.
\medskip

 This simple assertion gives information about divisors of the 
repunits: in particular, if $n$ is prime and  $A_n =p_1p_2\dots p_s$ is the 
factorization of $A_n$ into  prime  factors,  then  $h(p_1)=h(p_2)= \dots
=h(p_s)=n$. For instance, $A_7  = 239\cdot 4649$ and $h(239)=h(4649)=7$.
\medskip

{\bf Lemma 6:} Let $B_n(a,b)$ be an absolute prime, $p$ be a prime such that  $n>p-1$.
Suppose that
$10$ is  a primitive root modulo $p$, and $a$ and $p$ are relatively
prime. Then $n$  is a multiple of $p-1$.

{\bf Proof:} Let us consider the integers
\[
 B_i  = a\cdot 10^{p-1}\cdot A_{n-p+1} + a\cdot A_{p-1} + (b-a)\cdot
10^i,\qquad  0\le i\le p-2, 
\]
obtained from $B_n(a,b)$ by a permutation of  the  last 
$p-1$  digits. The  powers $10^i$, $0\le i\le p-2$, yield all nonzero remainders on
dividing by $p$, so do the integers $(b-a)\cdot 10^i$, $0\le i\le p-2$, and hence all
the integers $B_i$ can be simultaneously prime only in the case when the integer $L =
a\cdot 10^{p-1}\cdot A_{n-p+1} + a\cdot A_{p-1}$  is divisible  by $p$. But then, since
$\hbox{GCD}(a\cdot 10^{p-1},p)=1$ and $A_{p-1}\equiv 0\ \mbox{(mod $p$)}$,  it  follows  that 
$A_{n-p+1}$ is  divisible by $p$ and by Lemma~5 $n$ is divisible by $p-1$.
\medskip

{\bf Lemma 7:} The integers  $B_n(a,b)$, $a\ne b$, are  not  absolutely 
prime for $7\le n\le 16$.

{\bf Proof:} If $a\ne 7$, it follows from Lemma 6, applied for $p=7$, that  we  need  to 
verify the integers $B_n(a,b)$ with $n=12$ only, whereas the case $a=7$ 
requires a little bit more work. Direct calculations (or the use of a computer) here seem 
to  be unavoidable. These calculations show that the integers $B_n(7,b)$ are either
multiples  of $3$, or else by a permutation of digits they can be converted into
multiples of $17$ or $19$.
\medskip

{\bf Theorem 2:}  Let $N$  be  an  absolute  prime,  different  from 
repunits, that contains $n>3$ digits in its decimal  representation. 
Then $n$ is a multiple of $11088$.

{\bf Proof:} According to the previous lemma we assume that $n>16$. 
Since $10$ is a primitive root modulo $17$, Lemma 6 yields that $n$ 
divides $16$ and hence $n\ge 32$.  We  can  repeat  this  argument  three 
times, using the primes $19$, $23$, $29$, to obtain that $n$ is a multiple  of 
$18$, $22$ and $28$, respectively. Therefore $n$ divides $\hbox{LCM}(16,18,22,28) 
= 11088$.
\medskip

     Richert [5] used in addition the primes $47$, $59$, $61$, $97$, $167$, 
$179$, $263$, $383$, $503$, $863$, $887$, $983$ to show that  the  number  $n$  of 
digits of the  absolute  prime  number $B_n(a,b)$  is  divisible  by 
$321,653,308,662,329,838,581,993,760$. He also  mentioned,  that  by 
using the tables of primes and their primitive roots up to $10^5$, it 
is possible to show that $n > 6\cdot 10^{175}$.
\medskip

     Let us discuss now what pairs $(a,b)$ can appear  in  a  decimal 
representation of an absolute prime $B_n(a,b)$ with $n>3$ (if it exists at 
all!).
\medskip

{\bf Theorem 3:} If for $n>3$ the integer $B_n(a,b)$ is an absolute 
prime, then $(a,b) \ne (9,7), (9,1), (1,7), (7,1), (3,9), (9,3)$.

{\bf Proof:}  Let us write down the following equality
\[
9A_n - 2\cdot 10^r  = 10^n - 1 - 2\cdot 10^r  = 10^n  + 1 - 2(10^r+1).
\]        
We know from Theorem 2 that $n$ must be even. Write $n = 2^m\cdot u$,  where 
$u$ is odd. Then for $r=2^m$ the integer $10^n +1$ is  divisible  by  $10^r +1$, 
and the integer $9A_n - 2\cdot 10^r$ is composite. But this integer  can  be 
obtained by a permutation of digits of $B_n(9,7)$.

     Furthermore,
\[
B_n(9,1) = 9A_n- 8 = 10^n - 9 = (10^{n/2} - 3)(10^{n/2}+ 3),
\]     
and this number is also composite.

     Finally, since $n$ by Theorem 2 is divisible by 3, the sums  of 
the digits of $B_n(1,7)$ and $B_n(7,1)$ are also divisible by $3$.  Hence 
these numbers are composite as well as $B_n(9,3)$ and $B_n(3,9)$.

\bigskip
 
\bigskip 

\centerline{REFERENCES}                           

\begin{enumerate}
\item M. Gardner, Mathematical games, Scientific American, June, 
1964, p. 118. 
\item  H.S.  Williams,  Some  prime  with  interesting   digits 
patterns, Math. Comp., 1978, Vol. 32, pp.1306-1310. 
\item S. Yates, Repunits  and  repetends.  -  Star  Publ.  Co., 
Boynton Beach, Florida, 1982. 
\item P. Ribenboim, Prime numbers records. A  new  Chapter  for 
the Guinness Book 
\item H.-E. Richert, On permutable primtall,  Norsk  Matematiske 
Tiddskrift, 1951, vol.33, pp.50-54.
\item T.N. Bhargava and P.H. Doyle, On the existence of absolute 
primes, Mathematics Magazine, 1974, vol.47, p.233.
\item A.W. Johnson, Absolute primes, Mathematics Magazine, 1977, 
vol.50, pp. 100-103.
\end{enumerate}

\end{document}